\documentclass[runningheads]{llncs}
\usepackage[T1]{fontenc}
\usepackage{graphicx}
\usepackage{amsfonts}
\usepackage{amsmath}
\usepackage{soul,color}

\newcommand{\norm}[1]{\left\lVert#1\right\rVert}
\begin{document}
\title{Oscillatory behaviour of the RBF-FD approximation accuracy under increasing stencil size}
\titlerunning{Oscillations in the RBF-FD error dependence on stencil size}
\author{Andrej Kolar-Požun\inst{1,2}\orcidID{0000-0001-9694-3828} \and
	Mitja Jančič\inst{1,3}\orcidID{0000-0003-1850-412X} \and
	Miha Rot\inst{1,3}\orcidID{0000-0002-8869-4507} \and
	Gregor Kosec\inst{1}\orcidID{0000-0002-6381-9078}}
\authorrunning{A. Kolar-Požun et al.}
\institute{
	Jožef Stefan Institute, Parallel and Distributed Systems Laboratory, Ljubljana, Slovenia
		\email{\{andrej.pozun,mitja.jancic,miha.rot,gregor.kosec\}@ijs.si}  \and
	Faculty of Mathematics and Physics, University of Ljubljana, Slovenia \and
	Jožef Stefan International Postgraduate School, Ljubljana, Slovenia
}

\maketitle
\begin{abstract}
When solving partial differential equations on scattered nodes using the Radial Basis Function generated Finite Difference (RBF-FD) method, one of the parameters that must be chosen is the stencil size. Focusing on Polyharmonic Spline RBFs with monomial augmentation, we observe that it affects the approximation accuracy in a particularly interesting way - the solution error oscillates under increasing stencil size. We find that we can connect this behaviour with the spatial dependence of the signed approximation error. Based on this observation we are then able to introduce a numerical quantity that indicates whether a given stencil size is locally optimal.
\end{abstract}

\section{Introduction}
Radial Basis Function generated Finite Differences (RBF-FD) is a method for solving Partial Differential Equations (PDEs) on scattered nodes that has recently been increasing in popularity. It uses Radial Basis Functions (RBFs) to locally approximate a linear differential operator in a chosen neighbourhood, generalising the well known finite difference methods. This neighbourhood used for the approximation is referred to as the stencil of a given point and is commonly chosen to simply consist of its $n$ closest neighbours.

Among the different possible choices of a RBF used, the Polyharmonic Splines (PHS) with appropriate polynomial augmentation stand out due to the fact that they possess no shape parameter, eliminating all the hassle that comes with having to find its optimal value. PHS RBF-FD has been studied extensively and proved to work well in several different contexts \cite{use1,use2,PHS2,InsightPHS,Role2}. 
Unlike in the case of RBFs with a shape parameter, where the approximation order is determined by the stencil size \cite{RBFFormulas}, in PHS RBF-FD it is determined by the degree of the monomials included in the augmentation \cite{InsightPHS}. Despite that, the choice of an appropriate stencil size can have a substantional impact on the accuracy. More precisely, the accuracy of the method displays an oscillatory behaviour under increasing stencil size. 

In the remainder of the paper, we present this observation and our findings.  
Ideally we would like to be able to predict the stencil sizes that correspond to the accuracy minima or at least provide some indicator on whether a given stencil size is near the minimum.

The following section describes our problem setup along with the numerical solution procedure and in section 3 our results are discussed.
\section{Problem setup}
Our analyses are performed on the case of the Poisson equation
\begin{equation} \label{eqn:poisson}
\nabla^2 u( \textbf{x}) = f(\textbf{x}),
\end{equation}
where the domain is a disc $\Omega = \{\textbf{x} \in \mathbb{R}^2 : \norm{x-(0.5,0.5)} \leq 0.5 \}$.
We choose the function $f(\textbf{x})$ such, that the problem given by Equation (\ref{eqn:poisson}) has a known analytic solution. 
Concretely, we choose
\begin{align}
u(x,y) =& \sin (\pi x) \sin (\pi y), \\
f(x,y) =& -2 \pi^2 \sin (\pi x) \sin (\pi y)
\end{align}
with the Dirichlet boundary conditions given by a restriction of $u(\textbf{x})$ to the boundary $\partial \Omega$.

We discretise the domain with the discretisation distance $h=0.01$, first discretising the boundary and then the interior using the algorithm proposed in \cite{fill}.
The Laplacian is then discretised using the RBF-FD algorithm as described in \cite{explainRBF}, where we choose the radial cubics as our PHS ($\phi(r) = r^3$) augmented with monomials up to degree $m=3$, inclusive. This requires us to associate to each discretisation point $\textbf{x}_i$ its stencil, which we take to consist of its $n$ nearest neighbours.
We can now convert the PDE (\ref{eqn:poisson}) into a sparse linear system, which we then solve to obtain an approximate solution $\hat{u}(\textbf{x})$. The source code is readily available in our git repository\footnote{\url{https://gitlab.com/e62Lab/public/2023\_cp\_iccs\_stencil\_size\_ef{}fect}}.

\begin{figure}
	\centering
	\includegraphics[width=.5\textwidth]{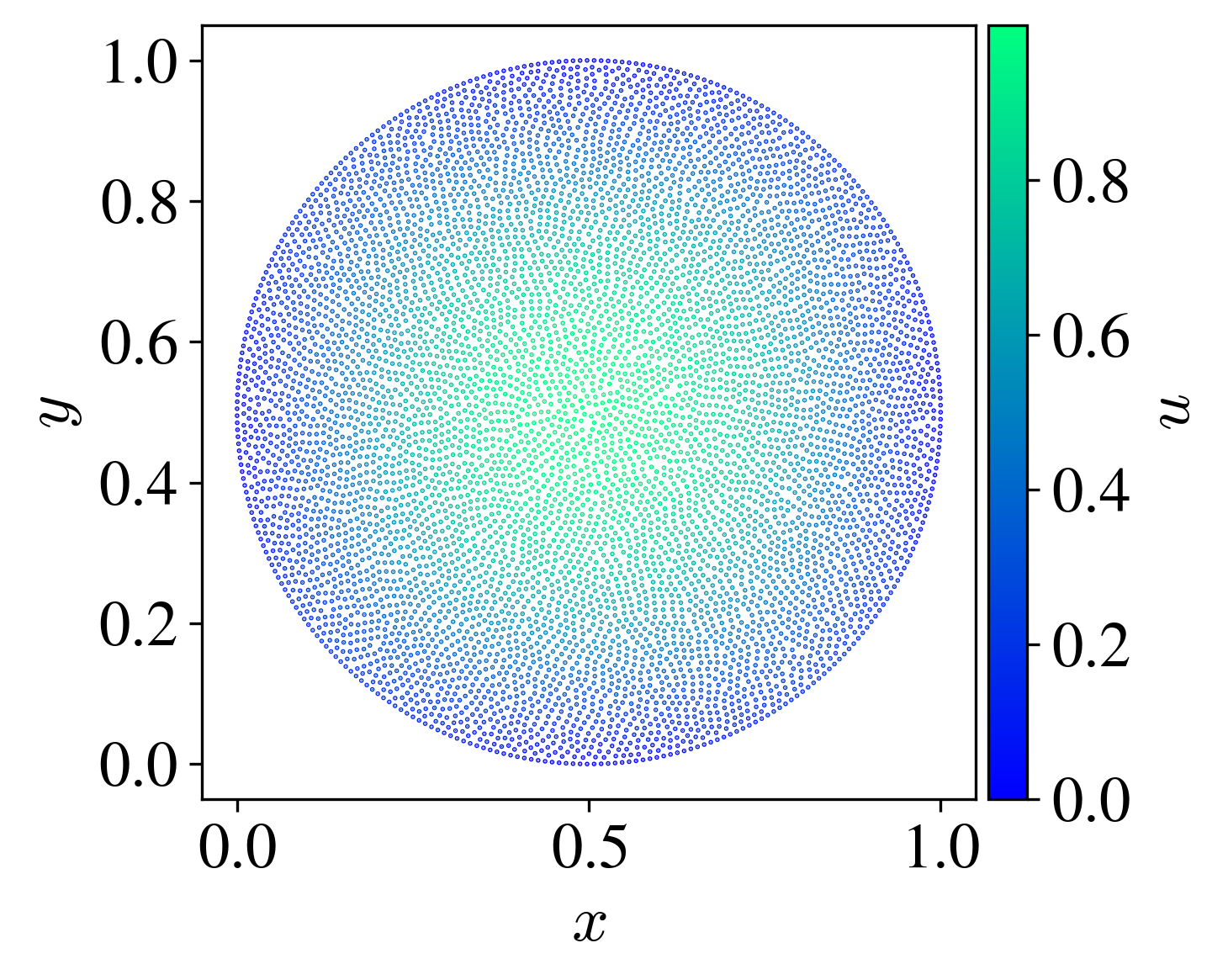}
	\caption{The analytical solution to the considered Poisson problem.}
	\label{fig:Solution1}
\end{figure}

The chosen analytical solution $u(\textbf{x})$ is displayed in Figure \ref{fig:Solution1}, which additionally serves as a visual representation of how the domain is discretised.

Having both the analytical and approximate solutions, we will be interested in the approximation error.
It will turn out to be useful to consider the signed pointwise errors of both the solution and the Laplacian approximation:
\begin{align}
&e_\mathrm{poiss}^\pm(\textbf{x}_i) = \hat{u}_i - u_i, \\
&e_\mathrm{lap}^\pm(\textbf{x}_i) = \tilde{\nabla}^2 u_i - f_i,
\end{align}
where $\tilde{\nabla}^2$ is the discrete approximation of the Laplacian and we have introduced the notation $u_i = u(\textbf{x}_i)$. The "poiss" and "lap" subscripts may be omitted in the text when referring to both errors at once.

As a quantitative measure of the approximation quality, we will also look at the average/max absolute value error:
\begin{align}
&e_\mathrm{poiss}^\mathrm{max} = \max_{\textbf{x}_i \in \mathring{\Omega}} |e_\mathrm{poiss}^\pm(\textbf{x}_i)|, \\
&e_\mathrm{poiss}^\mathrm{avg} = \frac{1}{N_\mathrm{int}} \sum_{\textbf{x}_i \in \mathring{\Omega}} |e_\mathrm{poiss}^\pm(\textbf{x}_i)|
\end{align}
and analogously for $e_\mathrm{lap}^\mathrm{max}$ and $e_\mathrm{lap}^\mathrm{avg}$. $N_\mathrm{int}$ is the number of discretisation points inside the domain interior $\mathring{\Omega}$.

In the next section we will calculate the approximation error for various stencil sizes $n$ and further investigate its (non-trivial) behaviour.

It is worth noting that the setup considered is almost as simple as it can be. The fact that we have decided not to consider a more complicated problem is intentional - there is no need to complicate the analysis by considering a more complex problem if the investigated phenomena already appears in a simpler one. This reinforces the idea that such behaviour arises from the properties of the methods used and not from the complexity of the problem itself.

\section{Results}
\begin{figure}
	\centering
	\includegraphics[width=\textwidth]{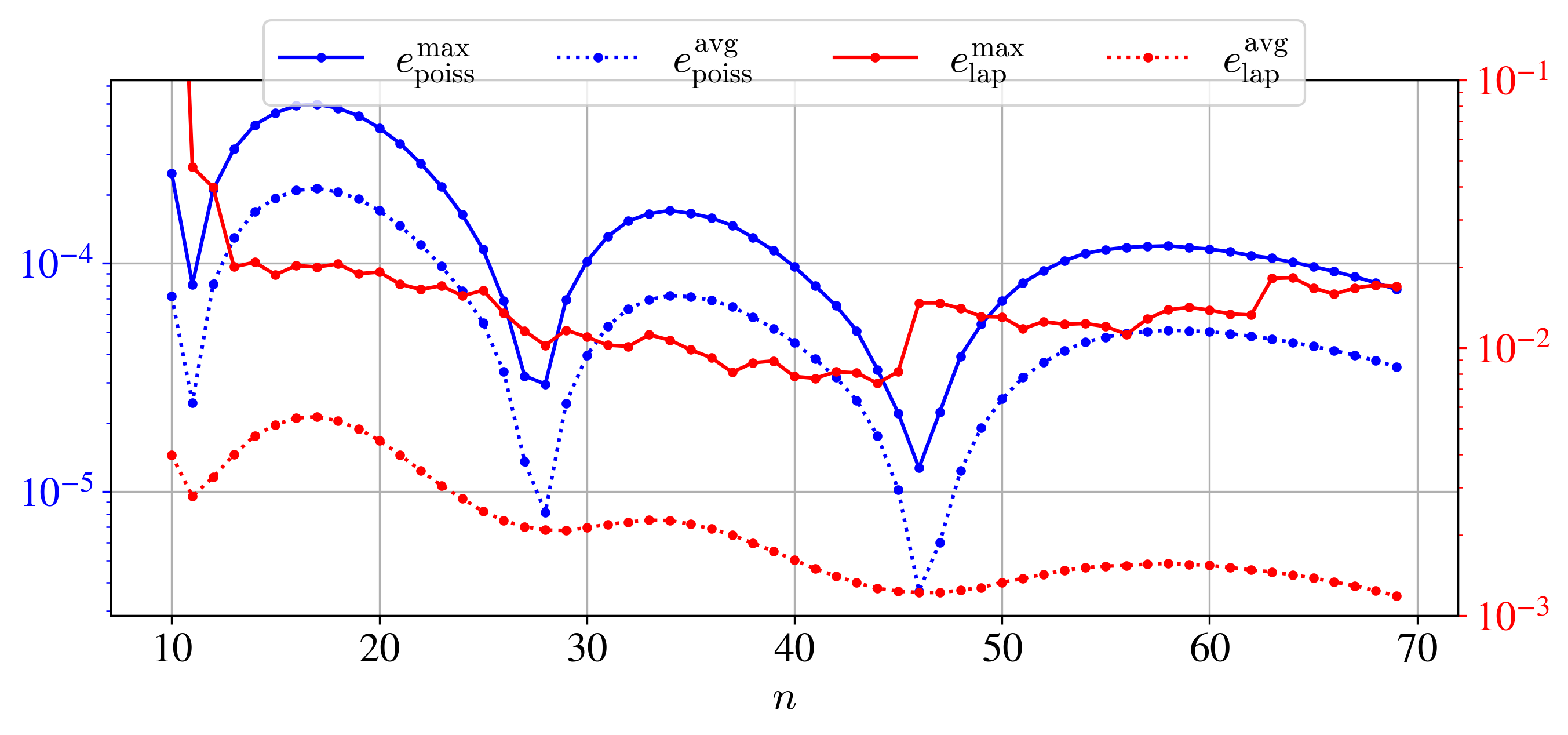}
	\caption{Dependence of the approximation errors on the stencil size $n$.}
	\label{fig:MeanErrors}
\end{figure}

In Figure~\ref{fig:MeanErrors} we see that $e_\mathrm{poiss}^\mathrm{max} (n)$ oscillates with several local minima (at stencil sizes $n=28,46$) and maxima (stencil sizes $n=17, 46$).
The dependence $e_\mathrm{poiss}^\mathrm{max} (n)$ seems to mostly resemble a smooth function. This is even more evident in $e_\mathrm{poiss}^\mathrm{avg}(n)$. The errors of the Laplacian are also plotted and we can observe that $e_\mathrm{lap}^\mathrm{avg}(n)$ has local minima and maxima at same stencil sizes. 
Such regularity implies that the existence of the minima in the error is not merely a coincidence, but a consequence of a certain mechanism that could be explained. Further understanding of this mechanism would be beneficial, as it could potentially allow us to predict the location of these local minima a priori. Considering that the error difference between the neighbouring local maxima and minima can be over an order of magnitude apart this could greatly increase the accuracy of the method without having to increase the order of the augmentation or the discretisation density. Note that the behaviour of $e_\mathrm{lap}^\mathrm{max}$ stands out as it is much more irregular. This implies that in order to explain the observed oscillations, we have to consider the collective behaviour of multiple points. This will be confirmed later on, when we consider the error's spatial dependence.

\begin{figure}
	\centering
	\includegraphics[width=\textwidth]{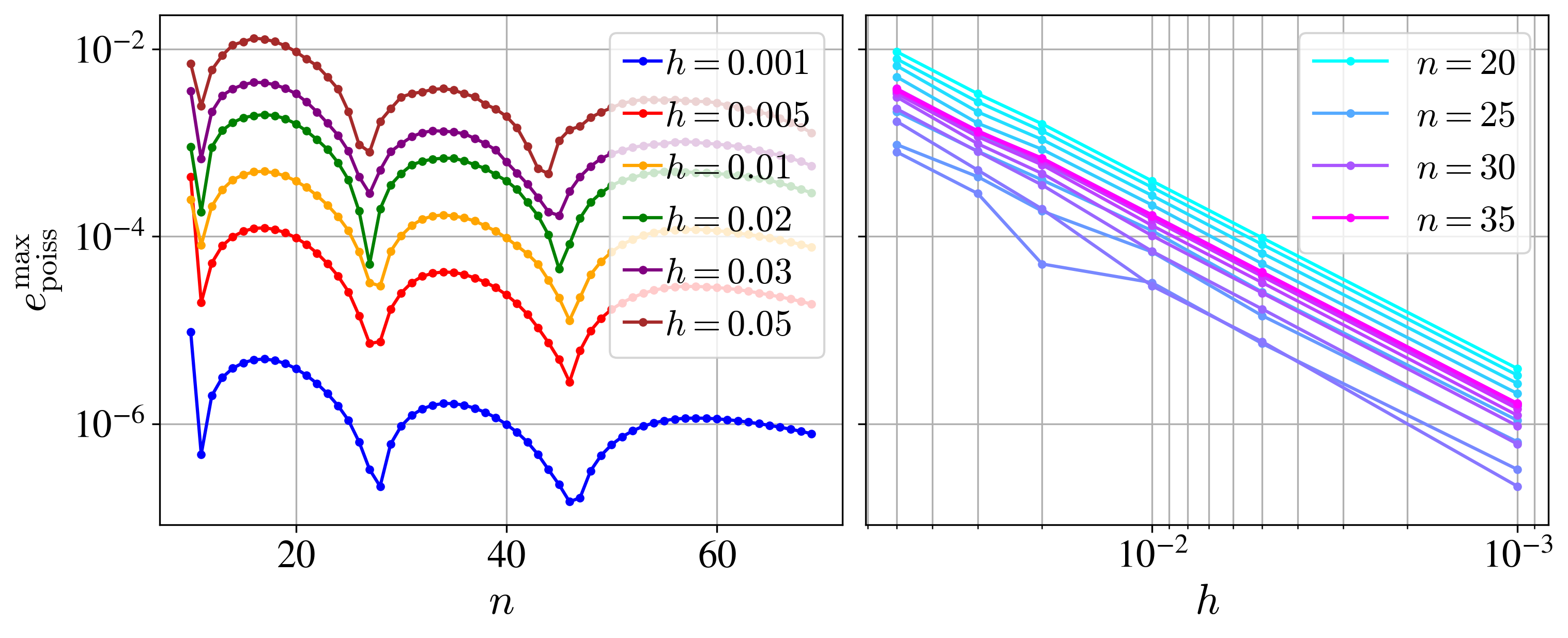}
	\caption{Behaviour of the approximation errors under a refinement of the discretisation.}
	\label{fig:Convergence}
\end{figure}

An immediate idea is that the choice of a sparse solver employed at the end of the solution procedure is responsible for the observed behaviour. We have eliminated this possibility by repeating the analysis with both the SparseLU and BiCGSTAB solvers, where no difference has been observed. The next idea we explore is the possibility of the discretisation being too coarse.
Figure~\ref{fig:Convergence} shows that under discretisation refinement $e_\mathrm{poiss}^\mathrm{max}(n)$ maintains the same shape and is just shifted vertically towards a lower error. The latter shift is expected, as we are dealing with a convergent method, for which the solution error behaves as $e \propto h^p$ as $h \to 0$, where $p$ is the order of the method. We also show $e_\mathrm{poiss}^\mathrm{max}(h)$ in a log-log scale for some different stencil sizes. It can be seen that the slopes and therefore the orders $p$ generally do not change with the stencil size and that the observed oscillations mainly affect the proportionality constant in $e \propto h^p$. The stencil dependence of the error proportionality constant has already been observed in similar methods \cite{RBFFormulas,layout}.

\begin{figure}
	\centering
	\includegraphics[width=.95\textwidth]{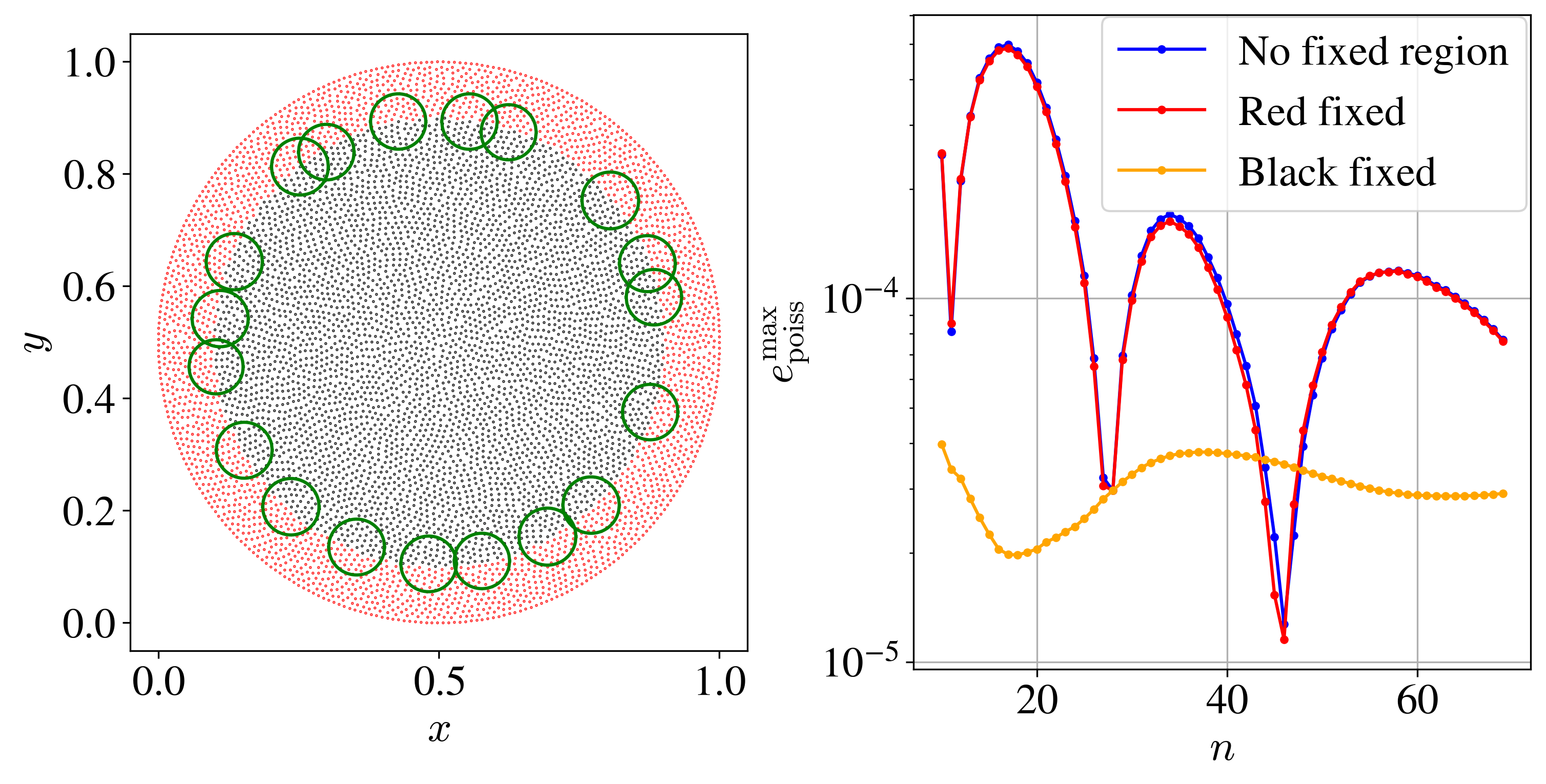}
	\caption{The seperation of the domain into two regions is seen on the left, where the green circles show the radii of the biggest stencils considered ($n=69$). The right graph shows the error dependence when either of the regions is at a fixed stencil size $n=28$. The previous result with no fixed stencil size regions is also shown.}
	\label{fig:BndFixed}
\end{figure}

Next we check if boundary stencils are responsible for the observed behaviour as it is known that they can be problematic due to their one-sideness \cite{Role2}. In Figure~\ref{fig:BndFixed} we have split our domain into two regions - the nodes near the boundary $\{\textbf{x}_i \in \Omega : \norm{\textbf{x}_i-(0.5,0.5)} > 0.4 \}$ are coloured red, while the nodes far from the boundary $\{\textbf{x}_i \in \Omega : \norm{\textbf{x}_i-(0.5,0.5)} \leq 0.4 \}$ are black.
We can see that the dependence of $e_\mathrm{poiss}^\mathrm{max} (n)$ marginally changes if we keep the stencil size near the boundary fixed at $n=28$ (corresponding to one of the previously mentioned minima), while only changing the stencil sizes of the remaining nodes. This shows that the observed phenomena is not a consequence of the particularly problematic boundary stencils.

\begin{figure}
	\centering
	\includegraphics[width=\textwidth]{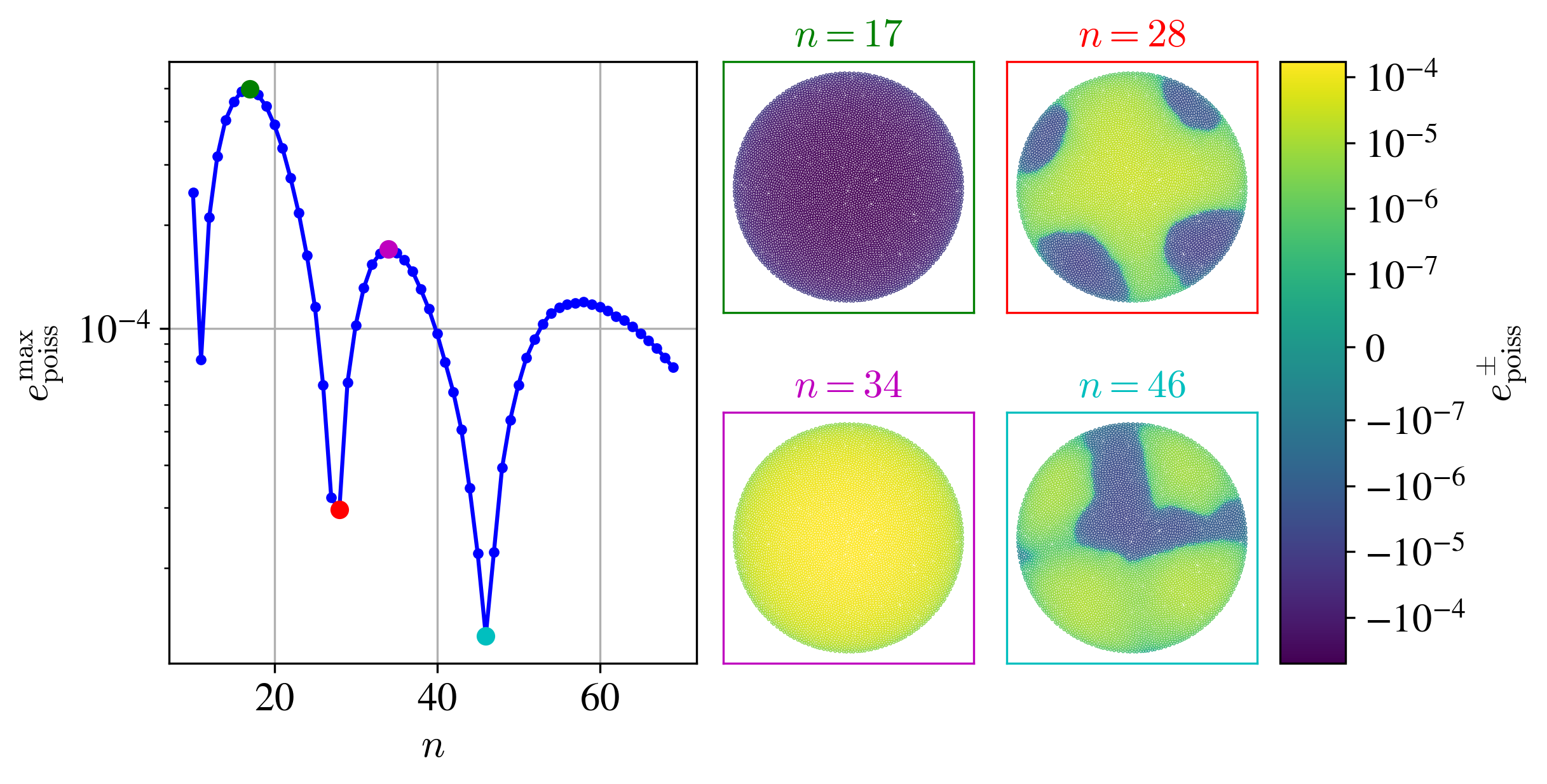}
	\caption{Spatial dependence of $e_\mathrm{poiss}^\pm$ in some local extrema. The colour scale is the same for all drawn plots.}
	\label{fig:PlotWithContours}
\end{figure}

Figure~\ref{fig:PlotWithContours} provides some more insight into the mechanism behind the oscillating error. Here we have plotted the spatial dependence of the signed error $e_\mathrm{poiss}^\pm$ for those stencils that correspond to the marked local extrema. We can observe that in the maxima, the error has the same sign throughout the whole domain. On the other hand, near the values of $n$ that correspond to the local minima there are parts of the domain that have differing error signs. Concretely, the sign of $e_\mathrm{poiss}^\pm$ is negative for stencil sizes between $17$ and $27$ inclusive. In the minima at $n=28$ both error signs are present, while for bigger stencil sizes (between $29$ and $45$ inclusive) the error again has constant sign only this time positive. The story turns out to repeat on the next minimum.

\begin{figure}
	\centering
	\includegraphics[width=\textwidth]{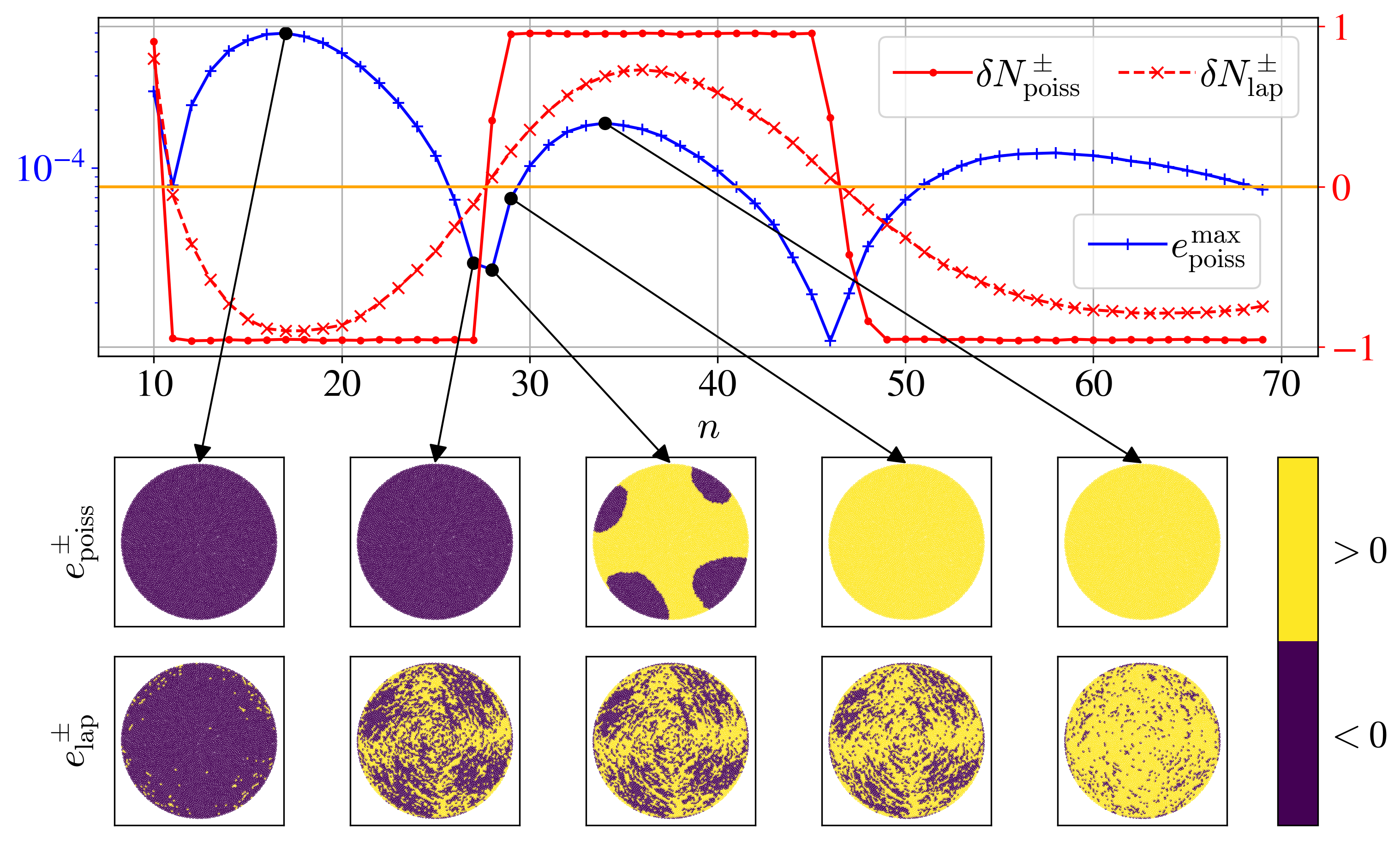}
	\caption{The quantities $\delta N^\pm (n)$ along with the spatial profiles of the signs of $e^\pm$ for some chosen stencil sizes. For convenience, $\delta N^\pm = 0$ is marked with an orange line.}
	\label{fig:MetricContours}
\end{figure}

This connection between the sign of $e_\mathrm{poiss}^\pm$ and the minima in $e_\mathrm{poiss}^\mathrm{max}(n)$ motivates us to define a new quantity:
\begin{equation} \label{eq:Metrika}
\delta N_\mathrm{poiss}^\pm = \frac{1}{N_\mathrm{int}}\left(|\{\textbf{x}_i \in \mathring{\Omega} : e_\mathrm{poiss}^\pm(\textbf{x}_i) > 0 \}| - |\{\textbf{x}_i \in \mathring{\Omega} : e_\mathrm{poiss}^\pm(\textbf{x}_i) < 0 \}| \right)
\end{equation}
and analogously for $\delta N_\mathrm{lap}^\pm$. Simply put, the quantity $\delta N_\mathrm{poiss}^\pm$ is proportional to the difference between the number of nodes with positively and negatively signed error. Assigning values of $\pm 1$ to positive/negative errors respectively, this quantity can be roughly interpreted as the average sign of the error. It should hold that $|\delta N_\mathrm{poiss}^\pm|$ is approximately equal to $1$ near the maxima and lowers in magnitude as we approach the minima. Figure~\ref{fig:MetricContours} confirms this intuition - $\delta N_\mathrm{poiss}^\pm (n)$ changes its values between $\pm 1$ very abruptly only near the $n$ that correspond to the minima of $e_\mathrm{poiss}^\mathrm{max} (n)$. A similar conclusion can be made for $\delta N_\mathrm{lap}^\pm$, which acts as a sort of "smoothed out" version of $\delta N_\mathrm{poiss}^\pm (n)$.

At a first glance, $N_\mathrm{lap}^\pm$ seems like a good candidate for an error indicator - it has a well-behaved dependence on $n$, approaches $\pm 1$ as we get closer to the error maxima and has a root near the error minima. The major downside that completely eliminates its applicability in the current state is the fact that we need access to the analytical solution to be able to compute it. 
\section{Conclusions}
Our study started with a simple observation - when solving a Poisson problem with PHS RBF-FD, the approximation accuracy depends on the stencil size $n$ in a non-trivial manner.
In particular, there exist certain stencil sizes where the method is especially accurate. A priori knowledge of these stencil sizes could decrease the solution error without any additional effort and is therefore strongly desirable. We have made a small step towards understanding this phenomena, eliminating various common numerical issues as the cause. 
Looking at the spatial dependence of the signed solution error, we have noticed that in the stencil sizes corresponding to the local error minima, the signed solution error is not strictly positive or negative. This is unlike the generic stencil sizes, where the error has the same sign throughout the domain. Motivated by this observation, we have introduced a quantity that is roughly the average sign of the pointwise Laplace operator errors and appears to have a root only near the stencils corresponding to the local error minima. 

The research presented is a step towards defining a more practically useful indicator, which would reveal the most accurate stencil sizes even without having access to the analytical solution and is a part of our ongoing research. Additional future work includes more rigorous theoretical explanations for the observations presented. Further experimental investigations should also be made, particularly to what extent our observations carry over to different problem setups - other differential equations and domain shapes.
\subsubsection{Acknowledgements}
The authors would like to acknowledge the financial support of Slovenian Research Agency (ARRS) in the framework of the research core funding No. P2-0095 and the Young Researcher programs PR-10468 and PR-12347.

\bibliographystyle{splncs04}
\bibliography{paper}

\end{document}